\documentclass[portuges,12pt,letter]{article}
\usepackage[centertags]{amsmath}
\usepackage{amsfonts}
\usepackage{newlfont}
\usepackage{amscd}
\usepackage{graphics}
\usepackage{epsfig}
\usepackage{indentfirst}
\usepackage{amsxtra}
\usepackage[latin1]{inputenc}
\usepackage{amssymb, amsmath}
\usepackage{amsthm}
\usepackage[mathscr]{eucal}

\newtheorem{thm}{Theorem}[section]

\title{\bf  A  convex dual formulation for a large class of non-convex models in variational optimization}
\date{}
\author{Fabio Silva Botelho \\  Department of Mathematics \\
Federal University of Santa Catarina \\
Florian\'{o}polis - SC, Brazil}
\begin{document}
\maketitle

\abstract{This short communication develops a  convex dual variational formulation for a large class of models in variational optimization. The  results are established through basic tools of functional analysis, convex analysis and duality theory. The main duality principle is developed as an application to a Ginzburg-Landau type system in superconductivity in the absence of a magnetic field. }

\section{Introduction}

In this section we establish a convex dual formulation for a large class of models in non-convex optimization.

The main duality principle is applied to the Ginzburg-Landau system in superconductivity in an absence of a magnetic field.

Such results are based on the works of J.J. Telega and W.R. Bielski \cite{2900,85,10,11} and on a D.C. optimization approach developed in Toland \cite{12}.

At this point we start to describe the primal and dual variational formulations.

Let $\Omega \subset \mathbb{R}^3$ be an open, bounded, connected set with a regular (Lipschitzian) boundary denoted by $\partial \Omega.$

For the primal formulation we consider the functional $J:U \rightarrow \mathbb{R}$ where
\begin{eqnarray}
J(u)&=& \frac{\gamma}{2}\int_\Omega \nabla u \cdot \nabla u\;dx \nonumber \\ && + \frac{\alpha}{2} \int_\Omega (u^2-\beta)^2\;dx -\langle u,f \rangle_{L^2}.
\end{eqnarray}

Here we assume $\alpha>0,\beta>0,\gamma>0$, $U=W_0^{1,2}(\Omega)$, $f \in L^2(\Omega)$. Moreover we denote
$$Y=Y^*=L^2(\Omega).$$

Define also $G:U  \rightarrow \mathbb{R}$ by
$$G(u)=\frac{\alpha}{2} \int_\Omega (u^2-\beta)^2\;dx+\frac{K}{2}\int_\Omega u^2\;dx,$$
and $F:U \rightarrow \mathbb{R}$ by

$$F(u)=-\frac{\gamma}{2}\int_\Omega \nabla u \cdot \nabla u \;dx+\frac{K}{2}\int_\Omega u^2\;dx+\langle u,f \rangle_{L^2},$$

It is worth highlighting that in such a case
$$J(u)=-F(u)+G(u), \; \forall u \in U.$$

From now and on, we assume a finite dimensional version for this model, in a finite elements of finite differences context, where,
for not relabeled operators and spaces, we also assume,
$$\gamma \nabla^2+ K > \mathbf{0}$$ in an appropriate matrices sense.

Furthermore, define  $$A^+=\{ u \in U\;:\; \delta^2J(u) \geq \mathbf{0}\}$$
$$(A^+)^0=\{ u \in U\;:\; \delta^2J(u) > \mathbf{0}\},$$
$$C^+=\{u \in U\;:\; uf \geq 0, \text{ in } \Omega\},$$
$$E^+=A^+ \cap C^+$$
and the following specific polar functionals specified, namely,
$G^*:Y^* \rightarrow \mathbb{R}$ by
\begin{eqnarray}G^*(v_1^*)&=&\sup_{ u \in A^+}\left\{\langle  u,  v_1^*\rangle_{L^2}-G(u)\right\}
\end{eqnarray}
and $F^*:Y^* \rightarrow \mathbb{R}$ by
\begin{eqnarray}F^*(v_1^*)&=&\sup_{ u \in U}\left\{ \langle  u, v_1^*\rangle_{L^2}-F(u)\right\}
\nonumber \\ &=& \frac{1}{2}\int_\Omega \frac{(v_1^*-f)^2}{\gamma \nabla^2+K}\;dx.
\end{eqnarray}

Define also
$J^*:Y^*  \rightarrow \mathbb{R}$ by

$$J^*(v_1^*)=F^*(v_1^*)-G^*(v_1^*)$$

Observe that there exists a Lagrange multiplier $\lambda \in W_0^{1,2}(\Omega)$ such that
\begin{eqnarray}
G^*(v_1^*)&=& \sup_{u \in U}\left\{ \langle u,v_1^*\rangle_{L^2}-G(u)+\frac{\gamma}{2}\int_\Omega \nabla \lambda \cdot \nabla \lambda\;dx\right. \nonumber \\ &&
\left.+\frac{6\alpha}{2}\int_\Omega \lambda^2u^2\;dx-\alpha\beta\int_\Omega \lambda^2\;dx\right\}.
\end{eqnarray}

Define now $G_2:Y^* \times U \times U \rightarrow \mathbb{R}$ by

$$G_2(v^*_1,u,\lambda)=\langle u,v_1^*\rangle_{L^2}-G(u)+\frac{\gamma}{2}\int_\Omega \nabla \lambda \cdot \nabla \lambda\;dx+\frac{6\alpha}{2}\int_\Omega \lambda^2u^2\;dx-\alpha\beta\int_\Omega \lambda^2\;dx.$$

Observe also that $$G^*(v_1^*)=G_2(v_1^*,\hat{u},\hat{\lambda}),$$

where $\hat{u}=u(v_1^*)$ and $\hat{\lambda}=\lambda(v_1^*)$ are such that

$$\frac{\partial G_2(v_1^*,\hat{u},\hat{\lambda})}{\partial u}=\mathbf{0},$$
and
$$\frac{\partial G_2(v_1^*,\hat{u},\hat{\lambda})}{\partial \lambda}=\mathbf{0}.$$

On the other hand,

$$\frac{\partial^2 G^*(v_1^*)}{\partial (v_1^*)^2}=\frac{\partial^2 G_2(v_1^*,\hat{u},\hat{\lambda})}{\partial (v_1^*)^2}
+\frac{\partial^2 G_2(v_1^*,\hat{u},\hat{\lambda})}{\partial v_1^*\partial u} \frac{\partial \hat{u}}{\partial v_1^*}+
+\frac{\partial^2 G_2(v_1^*,\hat{u},\hat{\lambda})}{\partial v_1^*\partial \lambda} \frac{\partial \hat{\lambda}}{\partial v_1^*}.$$

Moreover,

$$\frac{\partial^2 G_2(v_1^*,\hat{u},\hat{\lambda})}{\partial (v_1^*)^2}=\mathbf{0},$$
$$\frac{\partial^2 G_2(v_1^*,\hat{u},\hat{\lambda})}{\partial v_1^*\partial u}=1,$$
and
$$\frac{\partial^2 G_2(v_1^*,\hat{u},\hat{\lambda})}{\partial v_1^*\partial \lambda}=\mathbf{0}.$$

From these last results we get
$$\frac{\partial^2 G^*(v_1^*)}{\partial (v_1^*)^2}=\frac{\partial \hat{u}}{\partial v_1^*}.$$

However from $$\frac{\partial G_2(v_1^*,\hat{u},\hat{\lambda})}{\partial u}=\mathbf{0},$$
we have

$$v_1^*-2\alpha(\hat{u}^2-\beta)\hat{u}-K\hat{u}+6\alpha \hat{\lambda}^2 \hat{u}= \mathbf{0}$$

Taking the variation in $v_1^*$ in this last equation, we obtain

\begin{eqnarray}\label{US1110}&&1-6\alpha \hat{u}^2\frac{\partial \hat{u}}{\partial v_1^*} +2 \alpha \beta\frac{\partial \hat{u}}{\partial v_1^*}
\nonumber \\ &&-K\frac{\partial \hat{u}}{\partial v_1^*} +6\alpha \hat{\lambda}^2\frac{\partial \hat{u}}{\partial v_1^*}
+12\alpha \hat{\lambda}\frac{\partial \hat{\lambda}}{\partial v_1^*}\hat{u}=\mathbf{0}.\end{eqnarray}

On the other hand we must have also

$$\frac{\gamma}{2}\int_\Omega \nabla \hat{\lambda} \cdot \nabla \hat{\lambda}\;dx+\frac{1}{2}\int_\Omega 6 \alpha \hat{\lambda}^2\hat{u}^2\;dx-
\int_\Omega \alpha \beta \hat{\lambda}^2\;dx =0,$$ so that taking the variation in $v_1^*$ for this last equation and considering that
$$-\gamma \nabla^2 \hat{\lambda}+6 \alpha\hat{u}^2 \hat{\lambda}-2 \alpha \beta \hat{\lambda}=\mathbf{0},$$ we get

$$12 \alpha \hat{\lambda}^2\hat{u} \frac{\partial \hat{u}}{\partial v_1^*}=\mathbf{0}.$$

Hence if locally $ \hat{\lambda}^2\hat{u} \neq 0,$ then locally $$\frac{\partial \hat{u}}{\partial v_1^*}=0.$$

On the other hand if $ \hat{\lambda}^2\hat{u}= \mathbf{0},$ then from (\ref{US1110}) we have

$$\frac{\partial \hat{u}}{\partial v_1^*}=\frac{1}{6\alpha \hat{u}^2-2 \alpha \beta -6 \alpha \hat{\lambda}^2+K}.$$

Recalling that
$$\frac{\partial^2 G^*(v_1^*)}{\partial (v_1^*)^2}=\frac{\partial \hat{u}}{\partial v_1^*},$$

we have got
\begin{equation}\frac{\partial^2 G^*(v_1^*)}{\partial (v_1^*)^2}=\left\{\begin{array}{lr}
0, &\text { if } \hat{\lambda}^2\hat{u} \neq 0,
 \\
\frac{1}{6\alpha \hat{u}^2-2 \alpha \beta -6 \alpha \hat{\lambda}^2+K}, &  \text{ if }\hat{\lambda}^2\hat{u}=0.
 \end{array} \right.\end{equation}

 Observe also that
 $$\frac{\partial^2 J^*(v_1^*)}{\partial (v_1^*)^2}=\frac{\partial^2 F^*(v_1^*)}{\partial (v_1^*)^2}-\frac{\partial^2 G^*(v_1^*)}{\partial (v_1^*)^2}
=\frac{1}{\gamma\nabla^2+K}-\frac{\partial \hat{u}}{\partial v_1^*},$$
so that, for  $\hat{\lambda}^2\hat{u}=0$ we obtain
\begin{eqnarray}\frac{1}{\gamma\nabla^2+K}-\frac{\partial \hat{u}}{\partial v_1^*} &=&
\frac{1}{\gamma\nabla^2+K}-\frac{1}{6\alpha \hat{u}^2-2 \alpha \beta -6 \alpha \hat{\lambda}^2+K} \nonumber \\ &=&
\frac{-\gamma \nabla^2 -K +6\alpha \hat{u}^2-2 \alpha \beta -6 \alpha \hat{\lambda}^2+K}{(\gamma\nabla^2+K)(6\alpha \hat{u}^2-2 \alpha \beta -6 \alpha \hat{\lambda}^2+K)} \nonumber \\ &=& \frac{ \delta^2J(\hat{u}) -6 \alpha \hat{\lambda}^2}{(\gamma \nabla^2+K)(6\alpha \hat{u}^2-2 \alpha \beta -6 \alpha \hat{\lambda}^2+K)} \nonumber \\ &\geq& 0.
\end{eqnarray}

Summarizing,
\begin{equation}\frac{\partial^2 J^*(v_1^*)}{\partial (v_1^*)^2}=\left\{\begin{array}{lr}
\frac{1}{\gamma \nabla^2+K}, &\text { if } \hat{\lambda}^2\hat{u} \neq 0,
 \\
\frac{ \delta^2J(\hat{u}) -6 \alpha \hat{\lambda}^2}{(\gamma \nabla^2+K)(6\alpha \hat{u}^2-2 \alpha \beta -6 \alpha \hat{\lambda}^2+K)}, & \text{ if } \hat{\lambda}^2\hat{u}=0.
 \end{array} \right.\end{equation}

Hence, in any case, we have obtained $$\frac{\partial^2 J^*(v_1^*)}{\partial (v_1^*)^2} \geq \mathbf{0}, \; \forall v_1^* \in Y^*$$ so that $J^*$
is convex in $Y^*$.

\section{The main duality principle, a convex dual variational formulation}

Our main result is summarized by the following theorem.

\begin{thm} Considering the definitions and statements in the last section,  suppose also $\hat{v}^* \in Y^*$ is such that $$\delta J^*(\hat{v}^*)=\mathbf{0}.$$
 Assume also $$u_0=\frac{\partial F^*(\hat{v}_1^*)}{\partial v_1^*} \in  E^+ \cap (A^+)^0.$$

  Under such hypotheses, we have  $$\delta J(u_0)=\mathbf{0},$$  and
 \begin{eqnarray}
 J(u_0)&=&\inf_{u \in E^+}\left\{J(u)\right\} \nonumber \\ &=&
 \inf_{v_1^* \in Y^*}J^*(v_1^*)  \nonumber \\ &=&
 J_1^*(\hat{v}_1^*).
 \end{eqnarray}
\end{thm}
\begin{proof} From the hypothesis $$\frac{\partial J^*(\hat{v}_1^*)}{\partial v_1^*}=\mathbf{0}.$$  so that
\begin{eqnarray}\frac{\partial J^*(\hat{v}_1^*)}{\partial v_1^*}&=&\frac{\partial F^*(\hat{v}_1^*)}{\partial v_1^*}
-\frac{\partial G^*_1(\hat{v}_1^*)}{\partial v_1^*} =\mathbf{0}.
\end{eqnarray}

Since from the previous section we have got that $J^*$ is convex on $Y^*$, we may infer that
$$J^*(\hat{v}_1^*)=\inf_{v_1^* \in Y^*} J^*(v_1^*).$$

Also, from these last results,
 $$u_0-\frac{\partial G^*(\hat{v}_1^*)}{\partial v_1^*}=\mathbf{0},$$
 so that, since the restriction is not active in a neighborhood of $u_0$, from the Legendre transform properties, we obtain
 $$\hat{v}_1^*= \frac{\partial G(u_0)}{\partial u},$$
 and
$$ \hat{v}_1^*= \frac{\partial F(u_0)}{\partial u},$$
and thus
$$\mathbf{0}=\hat{v}_1^*-\hat{v}_1^*=-\frac{\partial F(u_0)}{\partial u}+\frac{\partial G(u_0)}{\partial u}=\delta J(u_0).$$

Summarizing $\delta J(u_0)=\mathbf{0}.$

Also from the Legendre transform properties we have

$$F^*(\hat{v}_1^*)=\langle u_0,\hat{v}_1^* \rangle_{L^2}-F(u_0),$$
and
$$G^*(\hat{v}_1^*)=\langle u_0,\hat{v}_1^* \rangle_{L^2}-G(u_0),$$
so that
$$J^*(\hat{v}_1^*)=F^*(\hat{v}_1^*)-G^*(\hat{v}_1^*)=-F(u_0)+G(u_0)=J(u_0).$$

Finally, from similar results in \cite{700}, we may infer that $E^+$ is convex so that from this and $\delta J(u_0)=\mathbf{0}$, we get
$$J(u_0)=\min_{u \in E^+} J(u).$$

Joining the pieces, we have got
 \begin{eqnarray}
 J(u_0)&=&\inf_{u \in E^+}\left\{J(u)\right\} \nonumber \\ &=&
 \inf_{v_1^* \in Y^*}J^*(v_1^*)  \nonumber \\ &=&
 J^*(\hat{v}_1^*).
 \end{eqnarray}

The proof is complete.
\end{proof}

\end{document}